\documentclass{amsart}
\usepackage{amsmath,amssymb}
\usepackage[dvips]{graphics}
\usepackage[all]{xy}
\newtheorem{Theorem}{Theorem}[section]

\newtheorem{Proposition}[Theorem]{Proposition}
\newtheorem{Corollary}[Theorem]{Corollary}

\newtheorem{Example}[Theorem]{Example}

\newtheorem{Remark}[Theorem]{Remark}

\makeatletter
\@addtoreset{figure}{section}
\def\@thmcountersep{-}
\makeatother

\numberwithin{equation}{section}

\begin{document}

\title{Plane curves in an immersed graph in ${\mathbb R}^2$}

\author{Marisa Sakamoto}
\address{Department of Mathematics, School of Education, Waseda University, Nishi-Waseda 1-6-1, Shinjuku-ku, Tokyo, 169-8050, Japan}
\email{marisa@fuji.waseda.jp}

\author{Kouki Taniyama}
\address{Department of Mathematics, School of Education, Waseda University, Nishi-Waseda 1-6-1, Shinjuku-ku, Tokyo, 169-8050, Japan}
\email{taniyama@waseda.jp}
\thanks{The second author was partially supported by Grant-in-Aid for Scientific Research (C) (No. 24540100), Japan Society for the Promotion of Science.}


\date{}

\dedicatory{}


\begin{abstract}
For any chord diagram on a circle there exists a complete graph on sufficiently many vertices such that any generic immersion of it to the plane  contains a plane closed curve whose chord diagram contains the given chord diagram as a sub-chord diagram. For any generic immersion of the complete graph on six vertices to the plane the sum of averaged invariants of all Hamiltonian plane curves in it is congruent to one quarter modulo one half. 
\end{abstract}

\maketitle

\section{Introduction} 

Throughout this paper we work in the piecewise linear category. 
By $K_n$ we denote the complete graph on $n$ vertices. It is shown in \cite{C-G} and \cite{Sachs} that every embedding of $K_6$ into the $3$-dimensional Euclidean space ${\mathbb R}^3$ always contains a non-splittable 2-component link. It is also shown in \cite{C-G} that every embedding of $K_7$ into ${\mathbb R}^3$ always contains a nontrivial knot. There are a number of subsequent results. See for example \cite{Negami}, \cite{K-S}, \cite{R-S-T}, \cite{F-P-F-N}, \cite{Foisy}, \cite{F-H-L-M}, \cite{O'donnol}, \cite{Nikkuni} and \cite{Nikkuni2}. These results show if an abstract graph $G$ is sufficiently complicated then any embedding of $G$ into ${\mathbb R}^3$ contains a complicated knot, link or knotted subgraph with prescribed properties. 

In this paper we show some analogous results in two-dimensions. Namely we consider a generic immersion of a finite graph into the Euclidean plane ${\mathbb R}^2$ and consider the immersed plane closed curves contained in the immersed graph. 

Let $G$ be a finite graph. We consider $G$ as a topological space in the usual way. A {\it cycle} is a subgraph of $G$ that is homeomorphic to a circle. An $n$-cycle of $G$ is a cycle of $G$ that contains exactly $n$ vertices. We denote the set of all $n$-cycles of $G$ by $\Gamma_n(G)$ and the set of all  cycles of $G$ by $\Gamma(G)$. Namely 
$\displaystyle{
\Gamma(G)=\cup_{n\in{\mathbb N}}\Gamma_n(G). 
}$ 
By a {\it generic immersion} we mean a continuous map from $G$ to ${\mathbb R}^2$ whose multiple points are only finitely many transversal double points of interior points of edges. These double points are called {\it crossing points} and the number of crossing points of a generic immersion $f:G\to{\mathbb R}^2$ is denoted by $c(f)$. We consider a circle as a graph. A {\it plane curve} is a generic immersion of a circle. An $n$-chord diagram is a paired $2n$ points on an oriented circle. An $n$-chord diagram is also called a {\it chord diagram}. Two $n$-chord diagrams ${\mathcal C}_1$ on an oriented circle $S_1$ and ${\mathcal C}_2$ on an oriented circle $S_2$ are equivalent if there is an orientation preserving homeomorphism from $S_1$ to $S_2$ that maps each pair of points of ${\mathcal C}_1$ to a pair of points of  ${\mathcal C}_2$. We consider chord diagrams up to this equivalence relation unless otherwise stated. A {\it sub-chord diagram} of a chord diagram ${\mathcal C}$ on an oriented circle $S$ is a chord diagram ${\mathcal D}$ on $S$ whose paired points are paired points of ${\mathcal C}$.

We describe an $n$-chord diagram by a circle with $n$-dotted chords. Each chord joins paired points of the chord diagram. See for example Figure \ref{chord diagrams}. 

Let $S$ be an oriented circle and $f:S\to{\mathbb R}^2$ a plane curve. By ${\mathcal C}(f)$ we denote the chord diagram on $S$ whose paired points correspond to the multiple points of $f$. Thus ${\mathcal C}(f)$ is a $c(f)$-chord diagram on $S$. 
Let $G$ be a finite graph and $f:G\to{\mathbb R}^2$ a generic immersion. Then for each cycle $\gamma$ of $G$ we have a restriction map $f|_\gamma:\gamma\to{\mathbb R}^2$. Then we have a chord diagram ${\mathcal C}(f|_\gamma)$ on the circle $\gamma$. 
By ${\mathbb S}^1$ we denote the unit circle with counterclockwise orientation. 

\vskip 3mm

\begin{Theorem}\label{chord} Let $n$ be a natural number with $n\geq2$. Let ${\mathcal C}$ be an $n$-chord diagram on ${\mathbb S}^1$. Let $f:K_{4n}\to {\mathbb R}^2$ be a generic immersion. Then there is a cycle $\gamma\in\Gamma_{4n}(K_{4n})$ and a sub-chord diagram ${\mathcal D}$ of ${\mathcal C}(f|_\gamma)$ such that ${\mathcal D}$ is equivalent to ${\mathcal C}$. 
\end{Theorem}

\vskip 3mm

An application of Theorem \ref{chord} to knot projections is shown in Section \ref{knot projections}. 

Let $\pi:{\mathbb R}^3\to{\mathbb R}^2$ be the natural projection defined by $\pi(x,y,z)=(x,y)$. 
A {\it knot} is a subspace of ${\mathbb R}^3$ that is homeomorphic to a circle. 
A knot $K$ in ${\mathbb R}^3$ is said to be {\it in regular position} if the restriction map $\pi|_K:K\to{\mathbb R}^2$ is a generic immersion from the circle $K$ to ${\mathbb R}^2$. 

Let $S$ be a circle and $f:S\to{\mathbb R}^2$ a generic immersion. Let $a_2(f)$ be the averaged invariant of the second coefficient of Conway polynomial of knots \cite{Polyak} \cite{L-W}. Namely $a_2(f)$ is the average of $a_2(K)$ where $K$ varies over all $2^{c(f)}$ knots in ${\mathbb R}^3$ in regular position with $\pi(K)=f(S)$. Here $a_2(K)$ denotes the second coefficient of the Conway polynomial of a knot $K$. 
Let $f:K_6\to{\mathbb R}^2$ be a generic immersion. We define $\alpha(f)$ by 
\[
\alpha(f)=\sum_{\gamma\in\Gamma_6(K_6)}a_2(f|_\gamma).
\]

\vskip 3mm

\begin{Theorem}\label{average}
Let $f:K_6\to{\mathbb R}^2$ be a generic immersion. Then 
\[
\displaystyle{\alpha(f)  \equiv \frac{1}{4} \pmod{\frac{1}{2}}}.
\] 
\end{Theorem}

\vskip 3mm

\begin{Corollary}\label{K6}
Let $f:K_6\to{\mathbb R}^2$ be a generic immersion. Then there is a cycle $\gamma\in\Gamma_6(K_6)$ and a knot $K$ in ${\mathbb R}^3$ in regular position with $a_2(K)\neq0$ such that $f(\gamma)=\pi(K)$. 
\end{Corollary}

\vskip 3mm

Note that Corollary \ref{K6} implies the result \cite[Theorem 3.4]{Foisy2} that every generic immersion of $K_6$ contains a projection of a nontrivial knot. 

\vskip 3mm

\begin{figure}[htbp]
      \begin{center}
\scalebox{0.4}{\includegraphics*{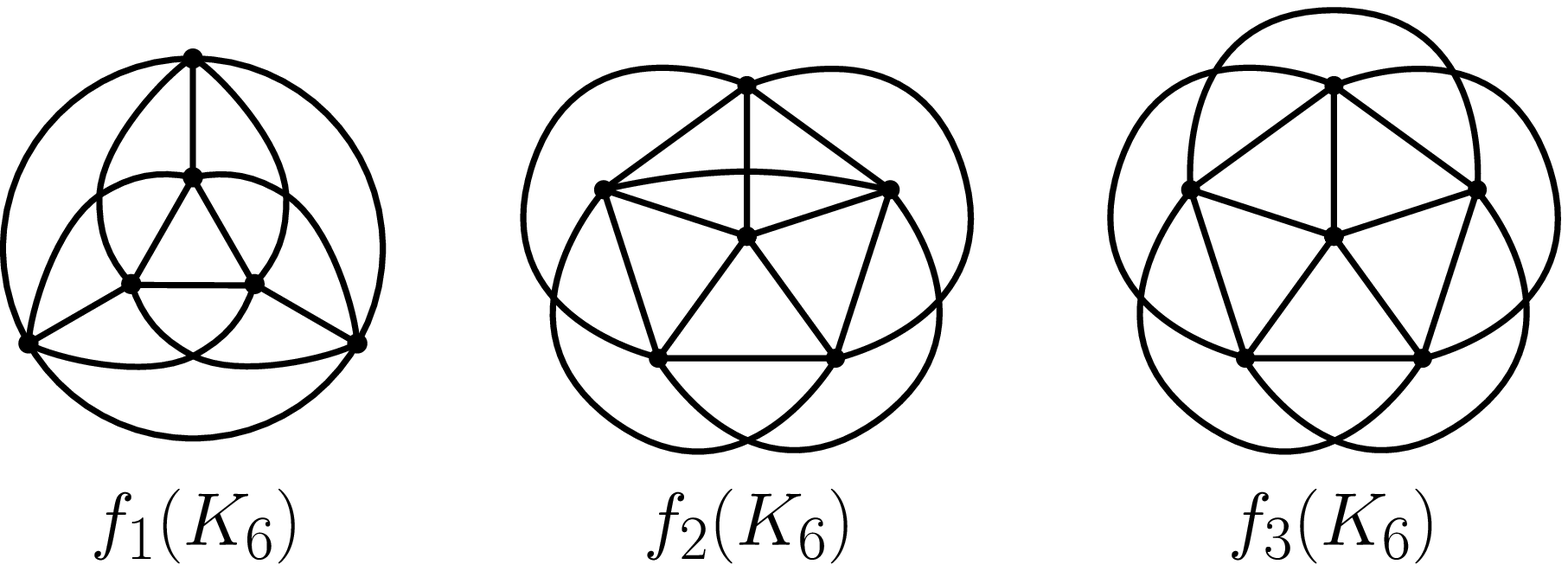}}
      \end{center}
   \caption{}
  \label{examples}
\end{figure}

\begin{Example}\label{example}
{\rm 
Let $f_1,f_2$ and $f_3$ be generic immersions from $K_6$ to ${\mathbb R}^2$ illustrated in Figure \ref{examples}. Then by a straightforward calculation we have $\alpha(f_1)=\frac{1}{4},\alpha(f_2)=\frac{3}{4}$ and $\alpha(f_3)=\frac{5}{4}$. 
}
\end{Example}

\section{Proofs}\label{proofs} 

We denote the set of all vertices of a graph $G$ by $V(G)$ and the set of all edges of $G$ by $E(G)$. 
Let $K_{m,n}$ be the complete bipartite graph on $m+n$ vertices partitioned into $m$ vertices and $n$ vertices. Let $G$ be the complete graph $K_5$ or the complete bipartite graph $K_{3,3}$. Let $f:G\to{\mathbb R}^2$ be a generic immersion. Let $x$ and $y$ be mutually disjoint edges of $G$. We denote the number of crossing points of $f$ contained in $f(x)\cap f(y)$ by $\#(f(x){\cap}f(y))$. We denote the set of all unordered pairs of mutually disjoint edges of $G$ by ${\mathcal E}(G)$. 
Set
\[
d(f)=\sum_{(x,y){\in}{\mathcal E}(G)}\#(f(x){\cap}f(y)).
\]
Namely $d(f)$ is the number of crossing points of $f$ made of disjoint edges of $G$. 

\vskip 3mm

\begin{Proposition}\label{K5K33} 
Let $G$ be the complete graph $K_5$ or the complete bipartite graph $K_{3,3}$. 
Let $f:G\to{\mathbb R}^2$ be a generic immersion. Then 
\[
d(f) \equiv 1 \pmod{2}.
\]
\end{Proposition}

\vskip 3mm

\begin{figure}[htbp]
      \begin{center}
\scalebox{0.35}{\includegraphics*{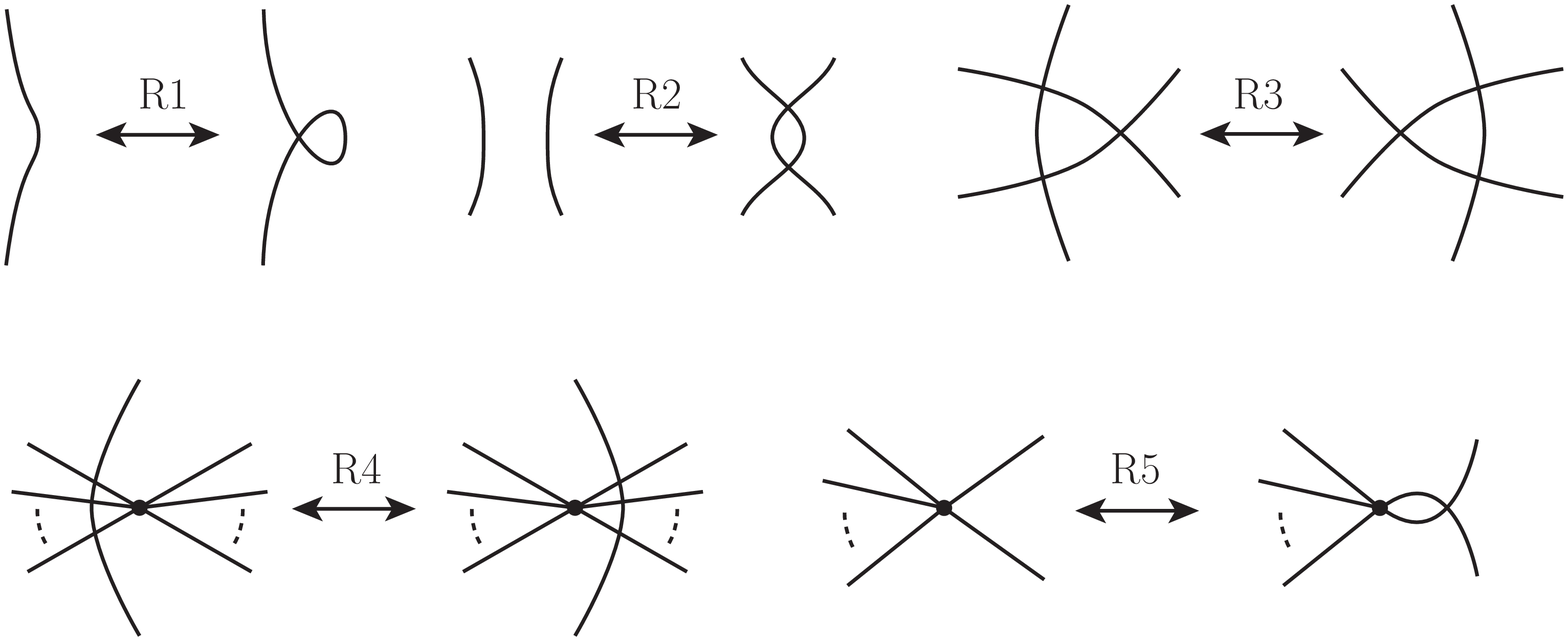}}
      \end{center}
   \caption{}
  \label{Reidemeister}
\end{figure}

\noindent{\bf Proof.} Let $f_0:G\to{\mathbb R}^2$ be a generic immersion illustrated in Figure \ref{K5K33}. Then we have 
\[
d(f_0)=1.
\]
It is well-known that any two generic immersions are transformed into each other up to self-homeomorphisms of $G$ and ${\mathbb R}^2$ by a finite sequence of the local moves illustrated in Figure \ref{Reidemeister}. It is easy to check that these moves do not change the parity of $d$ if $G=K_5$ or $G=K_{3,3}$. Thus we have the result. 
$\Box$

\begin{figure}[htbp]
      \begin{center}
\scalebox{0.4}{\includegraphics*{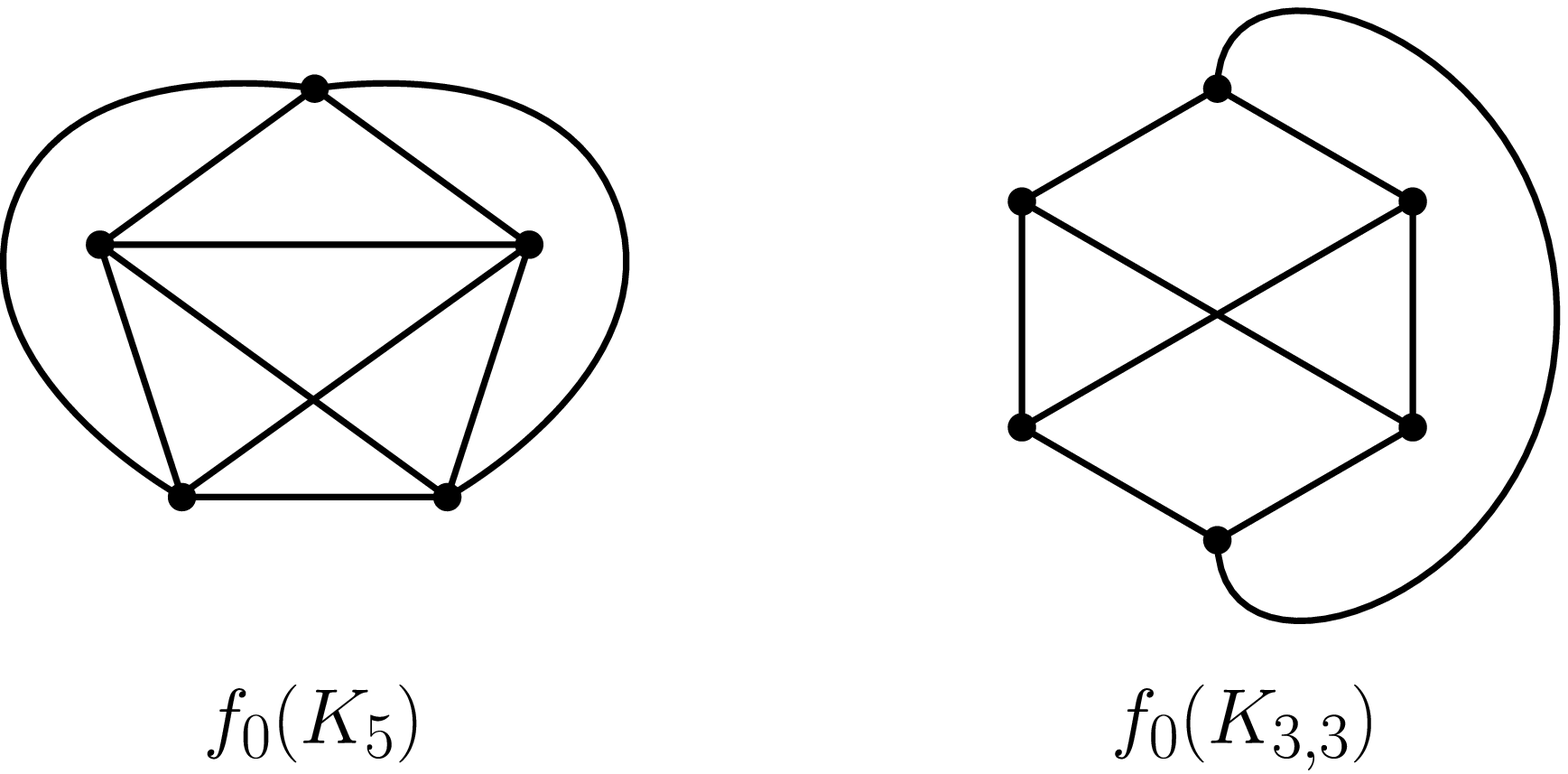}}
      \end{center}
   \caption{}
  \label{f0}
\end{figure} 

Let $G$ be a graph and $W$ a subset of $V(G)$. The induced subgraph $G[W]$ is the maximal subgraph of $G$ with $V(G[W])=W$. 

\vskip 3mm

\noindent{\bf Proof of Theorem \ref{chord}.} We will repeatedly apply Proposition \ref{K5K33} for certain subgraphs of $K_{4n}$ each of which is isomorphic to $K_5$ as follows. Set $V(K_{4n})=\{v_1,v_2,\cdots,v_{4n}\}$ and $H_1=K_{4n}[\{v_1,v_2,v_3,v_4,v_5\}]$. Then $H_1$ is isomorphic to $K_5$. By Proposition \ref{K5K33} there exist disjoint edges $x$ and $y$ of $H_1$ such that $f(x)\cap f(y)\neq\emptyset$. Let $c_1$ be a crossing in $f(x)\cap f(y)$. We may suppose without loss of generality that $x=v_1v_2$ and $y=v_3v_4$. Then set $H_2=K_{4n}[\{v_5,v_6,v_7,v_8,v_9\}]$. By the same argument we may suppose without loss of generality that $f(v_5v_6)\cap f(v_7v_8)\neq\emptyset$. Let $c_2$ be a crossing in $f(v_5v_6)\cap f(v_7v_8)$. Then set $H_3=K_{4n}[\{v_9,v_{10},v_{11},v_{12},v_{13}\}]$ and repeat the arguments. Thus we have $f(v_{4(k-1)+1}v_{4(k-1)+2})\cap f(v_{4(k-1)+3}v_{4(k-1)+4})\neq\emptyset$ for each $k\in\{1,2,\cdots,n-1\}$. Let $c_k$ be a crossing in $f(v_{4(k-1)+1}v_{4(k-1)+2})\cap f(v_{4(k-1)+3}v_{4(k-1)+4})$ for each $k\in\{1,2,\cdots,n-1\}$. Finally set the subgraph $H_n=K_{4n}[\{v_{4n-4},v_{4n-3},v_{4n-2},v_{4n-1},v_{4n}\}]$. There are two cases. Namely we may suppose without loss of generality that either $f(v_{4n-3}v_{4n-2})\cap f(v_{4n-1}v_{4n})\neq\emptyset$ or $f(v_{4n-4}v_{4n-3})\cap f(v_{4n-2}v_{4n-1})\neq\emptyset$. Let $c_n$ be a crossing in $f(v_{4n-3}v_{4n-2})\cap f(v_{4n-1}v_{4n})$ or $f(v_{4n-4}v_{4n-3})\cap f(v_{4n-2}v_{4n-1})$ respectively. 
In either case it is easy to find $\gamma\in\Gamma_{4n}(K_{4n})$ containing the edges $v_1v_2,v_3v_4,\cdots,v_{4n-7}v_{4n-6},v_{4n-5}v_{4n-4}$ and $v_{4n-3}v_{4n-2},v_{4n-1}v_{4n}$ or $v_{4n-4}v_{4n-3},v_{4n-2}v_{4n-1}$ such that the sub-chord diagram ${\mathcal D}$ of ${\mathcal C}(f|_\gamma)$ corresponding to the crossings $c_1,\cdots,c_n$ is equivalent to ${\mathcal C}$. This completes the proof. 
$\Box$

\vskip 3mm

\begin{figure}[htbp]
      \begin{center}
\scalebox{0.5}{\includegraphics*{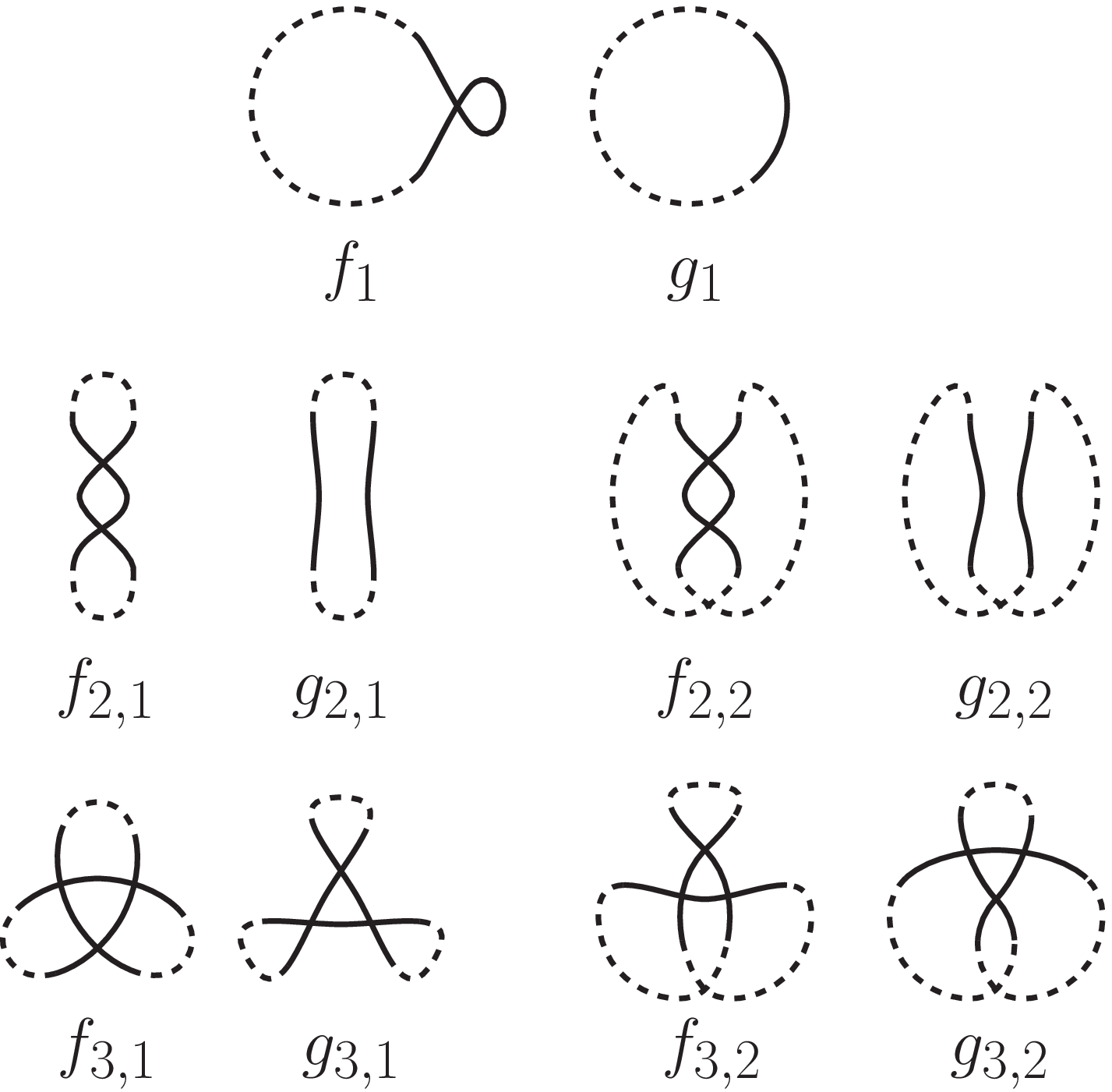}}
      \end{center}
   \caption{}
  \label{curves}
\end{figure}

\begin{Proposition}\label{average change} Let $f_1,g_1,f_{2,1},g_{2,1},f_{2,2},g_{2,2},f_{3,1},g_{3,1},f_{3,2}$ and $g_{3,2}$ be generic immersions from ${\mathbb S}^1$ to ${\mathbb R}^2$ such that each of $f_1$ and $g_1$, $f_{2,1}$ and $g_{2,1}$, $f_{2,2}$ and $g_{2,2}$, $f_{3,1}$ and $g_{3,1}$ and $f_{3,2}$ and $g_{3,2}$ differ locally as illustrated in Figure \ref{curves}. Then we have the following formulas. 

(1) $a_2(f_1)-a_2(g_1)=0$,

(2) $a_2(f_{2,1})-a_2(g_{2,1})=0$,

(3) $a_2(f_{2,2})-a_2(g_{2,2})=\frac{1}{4}$,

(4) $a_2(f_{3,1})-a_2(g_{3,1})=\frac{1}{4}$,

(5) $a_2(f_{3,2})-a_2(g_{3,2})=\frac{1}{4}$.
\end{Proposition}

\noindent{\bf Proof.} It is known in \cite{Polyak} that $a_2=\frac{1}{8}(J^++2St)$ where $J^+$ and $St$ are plane curve invariants defined by Arnold \cite{Arnold}. Then the formulas follow from the definitions of these invariants. We note that a direct calculation based on the well known formula $a_2(K_+)-a_2(K_-)={\rm lk}(L_0)$ in \cite{Kauffman} where $K_+$, $K_-$ and $L_0$ are knots and a $2$-component link forming a skein triple and ${\rm lk}$ denotes the linking number is also straightforward. 
$\Box$

\vskip 3mm

\noindent{\bf Proof of Theorem \ref{average}.} Let $f_1:K_6\to{\mathbb R}^2$ be a generic immersion in Example \ref{example}. Then $\alpha(f_1)=\frac{1}{4}$. Let $f:K_6\to{\mathbb R}^2$ be any generic immersion. Then $f$ and $f_1$ are transformed into each other up to self-homeomorphisms of $G$ and ${\mathbb R}^2$ by a finite sequence of the local moves ${\rm R}1,{\rm R}2,{\rm R}3,{\rm R}4,{\rm R}5$ illustrated in Figure \ref{Reidemeister}. Therefore it is sufficient to show that the parity of $4\alpha$ is invariant under these five local moves. Let $f:K_6\to{\mathbb R}^2$ and $g:K_6\to{\mathbb R}^2$ be generic immersions that differ by an application of a local move ${\rm R}$ where ${\rm R}$ is one of ${\rm R}1,{\rm R}2,{\rm R}3,{\rm R}4$ and ${\rm R}5$. We consider the following four cases. 

Case 1. ${\rm R}={\rm R}1$ or ${\rm R}={\rm R}5$. It follows from Proposition \ref{average change} (1) that $a_2(f|_\gamma)$ is invariant under ${\rm R}1$ and ${\rm R}5$ for each $6$-cycle $\gamma\in\Gamma_6(K_6)$. Therefore $\alpha(f)$ is invariant under ${\rm R}1$ and ${\rm R}5$.

Case 2. ${\rm R}={\rm R}2$. Let $e_1$ and $e_2$ be the edges of $K_6$ involved with ${\rm R}2$. 

Case 2.1. $e_1=e_2$.  In this case the $4!$ $6$-cycles of $K_6$ containing $e_1$ are involved. Note that all of them differ as $f_{2,1}$ and $g_{2,1}$ in Figure \ref{curves}, or all of them differ as $f_{2,2}$ and $g_{2,2}$ in Figure \ref{curves}. In the former case we apply Proposition \ref{average change} (2) and we have $\alpha(f)=\alpha(g)$. In the latter case we apply Proposition \ref{average change} (3) and we have $\alpha(f)-\alpha(g)=\pm\frac{4!}{4}$. 
Since $4!$ is even the parity of $4\alpha(f)$ and $4\alpha(g)$ coincide. 

Case 2.2. $e_1$ and $e_2$ are adjacent. In this case the $3!$ $6$-cycles containing both $e_1$ and $e_2$ are involved. Since $3!$ is even we have the conclusion by the same arguments as Case 2.1. 

Case 2.3. $e_1$ and $e_2$ are disjoint. In this case the $12$ $6$-cycles containing both $e_1$ and $e_2$ are involved. Exactly half of them are as $f_{2,1}$ and $g_{2,1}$ in Figure \ref{curves}, and the other half of them are as $f_{2,2}$ and $g_{2,2}$ in Figure \ref{curves}. Since the one half of $12$ is even, we have the conclusion. 

Case 3. ${\rm R}={\rm R}3$. Let $e_1$, $e_2$ and $e_3$ be the edges of $K_6$ involved with ${\rm R}3$. 

Case 3.1. $e_1=e_2=e_3$. In this case the $4!$ $6$-cycles of $K_6$ containing $e_1$ are involved. Note that all of them differ as $f_{3,1}$ and $g_{3,1}$ in Figure \ref{curves}, or all of them differ as $f_{3,2}$ and $g_{3,2}$ in Figure \ref{curves}. In the former case we apply Proposition \ref{average change} (4) and we have $\alpha(f)=\alpha(g)=\pm\frac{4!}{4}$. In the latter case we apply Proposition \ref{average change} (5) and we have $\alpha(f)-\alpha(g)=\pm\frac{4!}{4}$. 
Since $4!$ is even the parity of $4\alpha(f)$ and $4\alpha(g)$ coincide. 

Case 3.2. Exactly two of $e_1,e_2$ and $e_3$ are the same edge. We may suppose without loss of generality that $e_1=e_2$ and $e_3\neq e_1$. First suppose that $e_1$ and $e_3$ are adjacent. Then there are $3!$ $6$-cycles of $K_6$ containing $e_1$ and $e_3$ and by applying Proposition \ref{average change} (4) or (5), we have the conclusion. Next suppose that $e_1$ and $e_3$ are disjoint. Then there are $12$ $6$-cycles of $K_6$ containing both of $e_1$ and $e_3$. Exactly half of them are as $f_{3,1}$ and $g_{3,1}$ in Figure \ref{curves}, and the other half of them are as $f_{3,2}$ and $g_{3,2}$ in Figure \ref{curves}. Thus we have the conclusion. 

Case 3.3. $e_1\neq e_2\neq e_3\neq e_1$. We consider the following four cases.

Case 3.3.1. $e_1\cup e_2\cup e_3$ is homeomorphic to a circle. In this case no $6$-cycle of $K_6$ are involved and we have the result. 

Case 3.3.2. $e_1\cap e_2\cap e_3\neq\emptyset$. In this case no $6$-cycle of $K_6$ are involved and we have the result. 

Case 3.3.3. $e_1\cup e_2\cup e_3$ is homeomorphic to a closed interval. In this case the $2$ $6$-cycles containing all of $e_1,e_2$ and $e_3$ are involved and both of them are as $f_{3,1}$ and $g_{3,1}$ in Figure \ref{curves}, or both of them are as $f_{3,2}$ and $g_{3,2}$ in Figure \ref{curves}. Thus we have the conclusion. 

Case 3.3.4. $e_1$ and $e_2$ are adjacent and $e_3$ is disjoint from $e_1\cup e_2$. In this case the $4$ $6$-cycles containing all of $e_1,e_2$ and $e_3$ are involved and exactly two of them are as $f_{3,1}$ and $g_{3,1}$ in Figure \ref{curves}, and the other two of them are as $f_{3,2}$ and $g_{3,2}$ in Figure \ref{curves}. Thus we have the conclusion. 

Case 3.3.5. No two of $e_1,e_2$ and $e_3$ are adjacent. In this case the $8$ $6$-cycles containing all of $e_1,e_2$ and $e_3$ are involved and exactly four of them are as $f_{3,1}$ and $g_{3,1}$ in Figure \ref{curves}, and the other four of them are as $f_{3,2}$ and $g_{3,2}$ in Figure \ref{curves}. Thus we have the conclusion. 

\begin{figure}[htbp]
      \begin{center}
\scalebox{0.4}{\includegraphics*{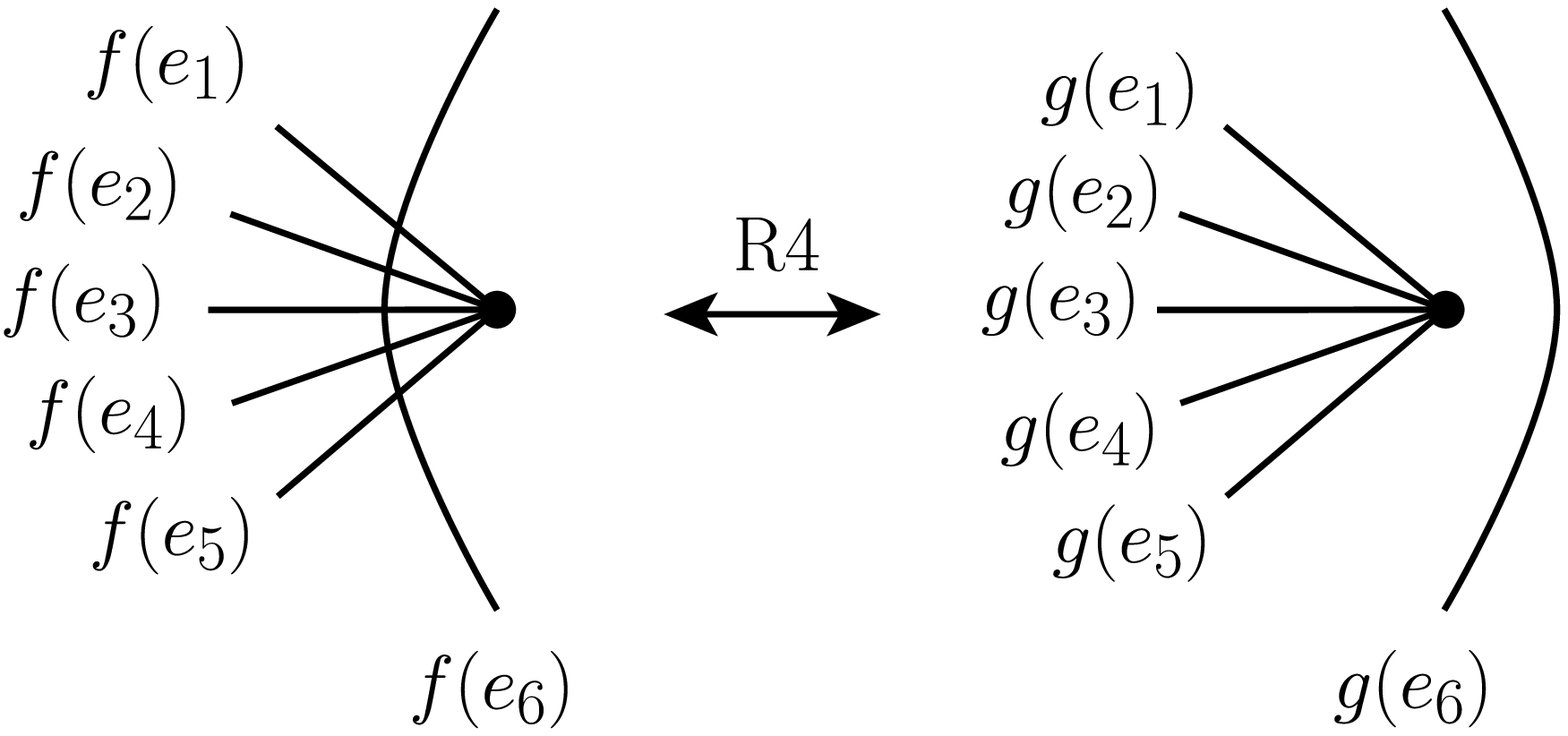}}
      \end{center}
   \caption{}
  \label{R4}
\end{figure}

Case 4. ${\rm R}={\rm R}4$. Up to local move ${\rm R}2$, we may suppose without loss of generality that the local move ${\rm R}4$ on $K_6$ is as illustrated in Figure \ref{R4}. Note that the edges $e_1,e_2,e_3,e_4$ and $e_5$ in Figure \ref{R4} are mutually distinct and $e_6$ may be equal to one of these five edges. For each pair $1\leq i<j\leq5$, let $\Gamma_{i,j}$ be the set of $6$-cycles of $K_6$ that contains all of the edges $e_i,e_j$ and $e_6$. Since 
\[
\alpha(f)-\alpha(g)=\sum_{\gamma\in\Gamma_6(K_6)}(a_2(f|_\gamma)-a_2(g|_\gamma))
\]
and $a_2(f|_\gamma)=a_2(g|_\gamma)$ for each $\gamma\in\Gamma_6(K_6)\setminus\bigcup_{1\leq i<j\leq5}\Gamma_{i,j}$, it is sufficient to show that 
\[
\sum_{\gamma\in\Gamma_{i,j}}4(a_2(f|_\gamma)-a_2(g|_\gamma))
\]
is even for each pair $1\leq i<j\leq5$. 

Case 4.1. $e_6=e_i$ or $e_6=e_j$. The $3!$ $6$-cycles containing both of $e_i$ and $e_j$ are involved. By applying Proposition \ref{average change} (2) or (3), we have the conclusion. 

Case 4.2. $e_6\neq e_i$ and $e_6\neq e_j$. If $e_i\cup e_j\cup e_6$ is homeomorphic to a circle or the complete bipartite graph $K_{1,3}$, then no $6$-cycle of $K_6$ is involved. Suppose that $e_i\cup e_j\cup e_6$ is homeomorphic to a closed interval. Then the two $6$-cycles of $K_6$ containing all of $e_i,e_j$ and $e_6$ are involved. By applying Proposition \ref{average change} (2) or (3), we have the conclusion. Suppose that $e_6$ is disjoint from $e_i\cup e_j$. Then the $4$ $6$-cycles containing all of $e_i,e_j$ and $e_6$ are involved. Exactly two of them are as $f_{2,1}$ and $g_{2,1}$ in Figure \ref{curves}, and the other two of them are as $f_{2,2}$ and $g_{2,2}$ in Figure \ref{curves}. Thus we have the conclusion. $\Box$

\vskip 3mm

\begin{figure}[htbp]
      \begin{center}
\scalebox{0.4}{\includegraphics*{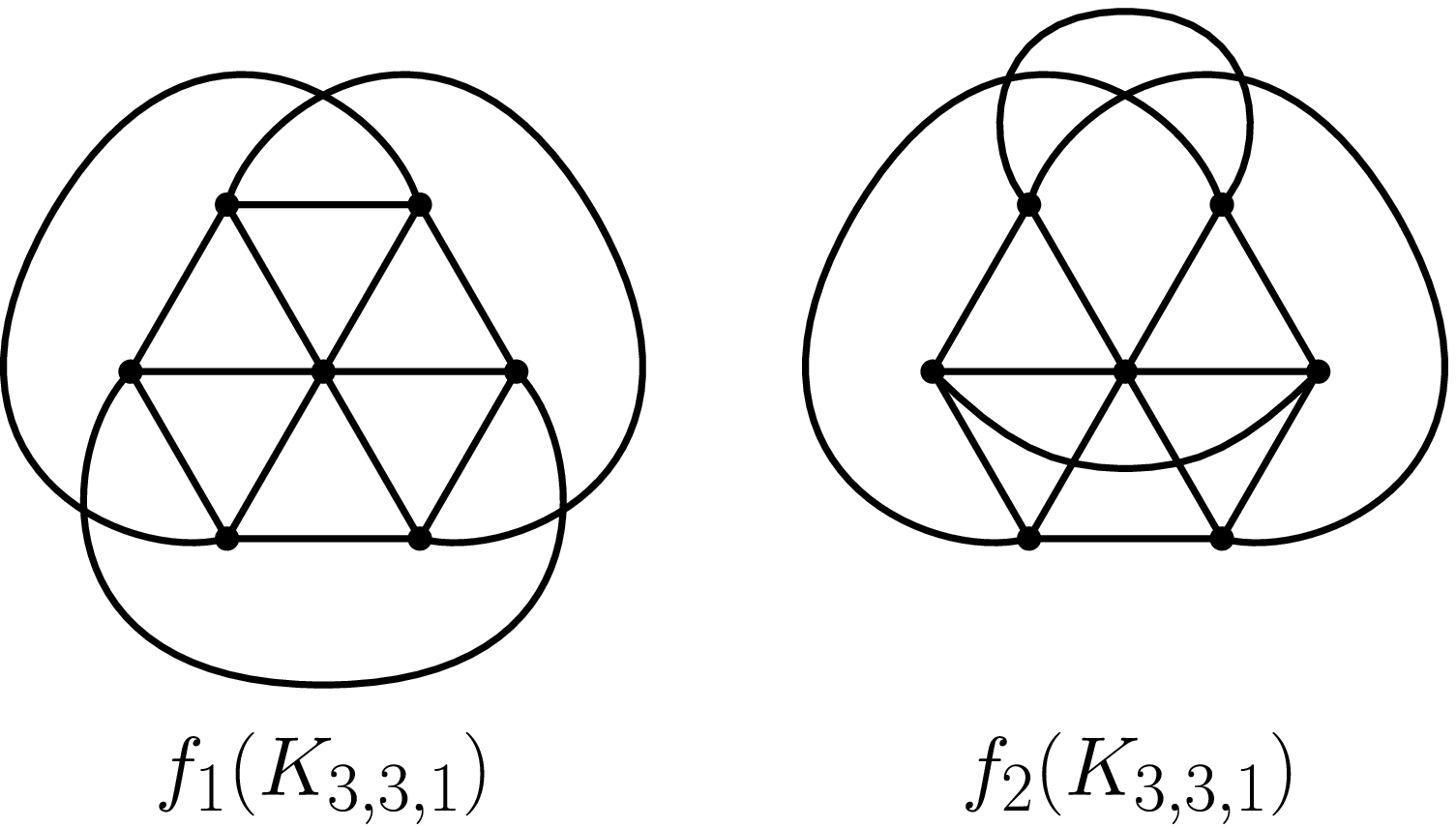}}
      \end{center}
   \caption{}
  \label{K331}
\end{figure}

\begin{Remark}\label{remark}
{\rm 
Let $K_{l,m,n}$ be the complete tripartite graph on $l+m+n$ vertices partitioned into $l$ vertices, $m$ vertices and $n$ vertices. 
It is known that as well as $K_6$, every embedding of the complete tripartite graph $K_{3,3,1}$ into the $3$-dimensional space ${\mathbb R}^3$ contains a non-splittable $2$-component link \cite{Sachs}. Here we show two generic immersions of $K_{3,3,1}$ to ${\mathbb R}^2$ that suggest non-existence of an analogy of Theorem \ref{average} for $K_{3,3,1}$. Let $f:K_{3,3,1}\to{\mathbb R}^2$ be a generic immersion. We define $\alpha(f)$ by the following. 
\[
\alpha(f)=\sum_{\gamma\in\Gamma_7(K_{3,3,1})}a_2(f|_\gamma).
\]
Let $f_1$ and $f_2$ be generic immersions of $K_{3,3,1}$ to ${\mathbb R}^2$ illustrated in Figure \ref{K331}. 
Then by a straightforward calculation we have $\alpha(f_1)=\frac{3}{4}$ and $\alpha(f_2)=1$. 
}
\end{Remark}

\section{Knot projections in a plane immersed graph}\label{knot projections} 

\begin{figure}[htbp]
      \begin{center}
\scalebox{0.4}{\includegraphics*{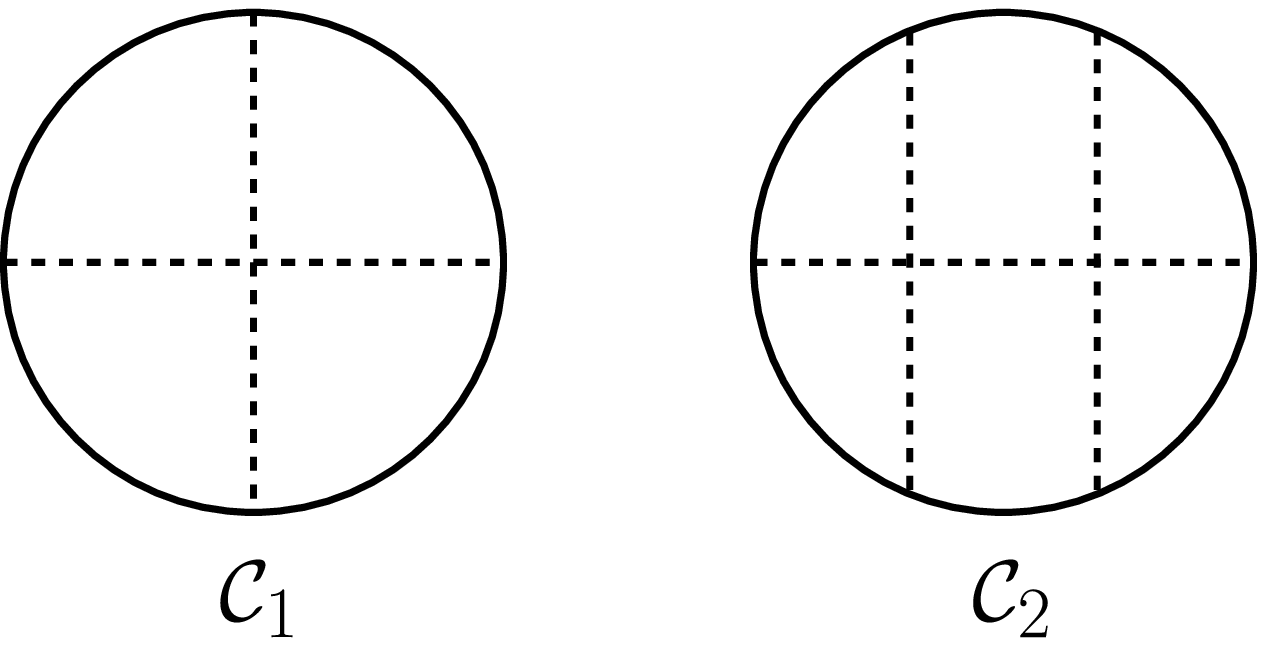}}
      \end{center}
   \caption{}
  \label{chord diagrams}
\end{figure}

In this section we show an application of Theorem \ref{chord} to knot projections. Two knots $K_1$ and $K_2$ are said to be {\it ambient isotopic} if there is an orientation preserving self-homeomorphism $h:{\mathbb R}^3\to{\mathbb R}^3$ such that $h(K_1)=K_2$. 
Let $K$ be a knot in ${\mathbb R}^3$. Let $f:{\mathbb S}^1\to{\mathbb R}^2$ a generic immersion. We say that $f$ is a {\it regular projection of $K$} if there exists a knot $K'$ in regular position that is ambient isotopic to $K$ such that $\pi(K')=f({\mathbb S}^1)$. Let ${\mathbb R}^2\cup\infty$ be the one-point compactification of ${\mathbb R}^2\cup\infty$ and $\iota:{\mathbb R}^2\to{\mathbb R}^2\cup\infty$ the natural inclusion map. Two plane curves $f:{\mathbb S}^1\to{\mathbb R}^2$ and $g:{\mathbb S}^1\to{\mathbb R}^2$ are said to be {\it spherically equivalent} if there exists an orientation preserving homeomorphism $\varphi:{\mathbb S}^1\to{\mathbb S}^1$ and an orientation preserving homeomorphism $h:{\mathbb R}^2\cup\infty\to{\mathbb R}^2\cup\infty$ such that $h\circ\iota\circ f=\iota\circ g\circ\varphi$. Let $f:{\mathbb S}^1\to{\mathbb R}^2$ be a plane curve and $C$ a circle in ${\mathbb R}^2$ transversely intersecting $f({\mathbb S}^1)$ at two points, say $P_1$ and $P_2$. Let $A$ be a simple arc in $C$ joining $P_1$ and $P_2$. Let $p:{\mathbb S}^1\to{\mathbb R}^2$ and $q:{\mathbb S}^1\to{\mathbb R}^2$ be plane curves such that $p({\mathbb S}^1)\cup q({\mathbb S}^1)=f({\mathbb S}^1)\cup A$, $p({\mathbb S}^1)\cap q({\mathbb S}^1)=A$, the orientation of $p$ coincides with that of $f$ on $p({\mathbb S}^1)\setminus A$ and the orientation of $q$ coincides with that of $f$ on $q({\mathbb S}^1)\setminus A$. Then we say that $f$ is a {\it connected sum of $p$ and $q$}. 
Let ${\mathcal C}_1$ and ${\mathcal C}_2$ be chord diagrams illustrated in Figure \ref{chord diagrams}. The following Theorem \ref{trefoil} and Theorem \ref{figure eight} are implicitly shown in the proofs of Theorem 1 and Theorem 2 of \cite{Taniyama} respectively. 

\vskip 3mm

\begin{Theorem}\label{trefoil}
Let $f:{\mathbb S}^1\to{\mathbb R}^2$ be a generic immersion. Then the following conditions are mutually equivalent. 

(1) $f$ is a regular projection of some nontrivial knot $K$.

(2) $f$ is a regular projection of the trefoil knot. 

(3) $f({\mathbb S}^1)$ is not obtained from an embedded circle in ${\mathbb R}^2$ by repeated applications of the local move ${\rm R}_1$ illustrated in Figure \ref{Reidemeister}.

(4) The chord diagram ${\mathcal C}(f)$ contains a sub-chord diagram that is equivalent to ${\mathcal C}_1$. 
\end{Theorem}

\vskip 3mm

\begin{figure}[htbp]
      \begin{center}
\scalebox{0.5}{\includegraphics*{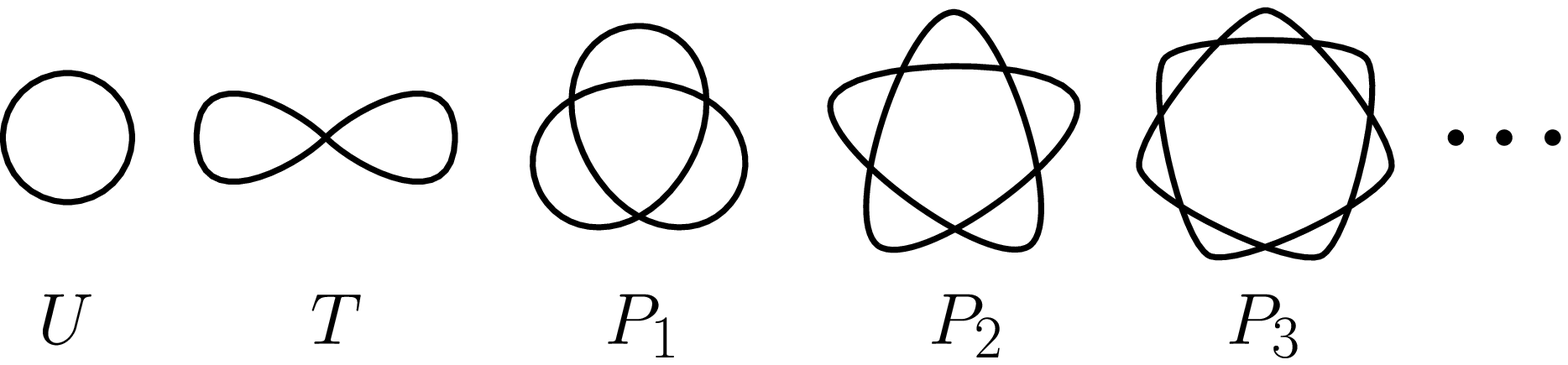}}
      \end{center}
   \caption{}
  \label{projections}
\end{figure}

\begin{Theorem}\label{figure eight}
Let $f:{\mathbb S}^1\to{\mathbb R}^2$ be a generic immersion. Then the following conditions are mutually equivalent. 

(1) $f$ is a regular projection of the figure-eight knot. 

(2) $f({\mathbb S}^1)$ is not equivalent to any connected sum of some plane curves each of which is equivalent to one of the plane curves $U,T,P_1,P_2,P_3,\cdots$ illustrated in Figure \ref{projections}. 

(3) The chord diagram ${\mathcal C}(f)$ contains a sub-chord diagram that is equivalent to ${\mathcal C}_2$. 
\end{Theorem}


Note that Theorem \ref{trefoil} and \cite[Theorem 3.4]{Foisy2} implies that for any generic immersion $f:K_6\to{\mathbb R}^2$ there exists a $6$-cycle $\gamma$ of $K_6$ such that $f(\gamma)$ is a regular projection of the trefoil knot. By Theorem \ref{chord} and Theorem \ref{figure eight} we have the following corollary.

\vskip 3mm

\begin{Corollary}\label{figure eight projection}
Let $f:K_{12}\to{\mathbb R}^2$ be a generic immersion. Then there is a cycle $\gamma\in\Gamma_{12}(K_{12})$ such that $f(\gamma)$ is a regular projection of the figure-eight knot. 
\end{Corollary}

\vskip 3mm


\section*{Acknowledgments} The authors are grateful to Professors Shin'ichi Suzuki and Toshie Takata for their encouragements. The authors are also grateful to Professors Toshiki Endo, Reiko Shinjo, Noboru Ito and Ryo Hanaki for their helpful comments.

{\normalsize
\end{document}